\documentclass{article}
\usepackage[normalem]{ulem}
\usepackage[most]{tcolorbox}
\usepackage{pdflscape}
\usepackage{hyperref}
\usepackage[
backend=biber,
giveninits=true, 
doi=false, 
url=false, 
isbn=false, 
eprint=true,
maxalphanames=99,
maxnames=99,
style=alphabetic,
citestyle=alphabetic]
{biblatex}
\DeclareFieldFormat{extraalpha}{#1}
\DeclareLabelalphaTemplate{
  \labelelement{
    \field[final]{shorthand}
    \field{label}
    \field[strwidth=3,strside=left,ifnames=1]{labelname}
    \field[strwidth=1,strside=left]{labelname}
  }
}
\usepackage{dsfont}
\usepackage{amsmath}
\usepackage{amsfonts}
\usepackage{amssymb}
\usepackage{enumerate}
\usepackage{amscd}
\usepackage{graphicx}
\usepackage{mathrsfs}
\usepackage{mathtools}
\usepackage{mathdots}

\usepackage{tikz}
\usepackage{ifthen}
\usepackage{tikz-cd}
\usepackage{quiver}
\usepackage{float}
\usetikzlibrary{decorations.markings}
\usetikzlibrary{arrows.meta}
\usetikzlibrary{positioning}
\usetikzlibrary{calc}
\usepackage{xfp}

\title{Arbitrary classes in $\geq 3$-degree cohomology of a finite group with arbitrary coefficients may be trivialized in a finite extension}
\author{Adrien DeLazzer Meunier}
\date{}

\def\Q{\mathbb{Q}}
\def\Z{\mathbb{Z}}

\newcommand{\bigslant}[2]{#1/#2}

\newcounter{Def}

\newenvironment{Definitions}{ \refstepcounter{Def}\mbox{}\\
\noindent\textbf{Definitions \arabic{Def}} }
{\vspace{0.3cm}}

\newenvironment{Definition}{ \refstepcounter{Def}\mbox{}\\
\noindent\textbf{Definition \arabic{Def}} }
{\vspace{0.3cm}}

\newenvironment{Proof}[1][\quad]
{\mbox{}\\\noindent \textbf{Proof:{\hspace{0.7cm}
#1}\\}}{\hfill$\Box$\\}

\newtheorem{Examples}[Def]{Examples}

\newtheorem{Remark}[Def]{Remark}
\newtheorem{Lemma}[Def]{Lemma}
\newtheorem{Theorem}[Def]{Theorem}
\newtheorem{Cor}[Def]{Corollary}

\newtheorem{Facts}[Def]{Facts}

\begin{document}

\maketitle

\section{Introduction}

The purpose of this note is to provide exposition for a proof of the statement in the title. This idea, that arbitrary cohomology classes (of high enough degree) of a finite group $G$ can be trivialized in a \textit{finite} group extension, has been known to experts for some time. In particular, this problem has been the subject of discussion in the physics literature as in, for example,\cite{YujiRecipe} which serves as the primary inspiration for this note. The author also points to \cite{WWW} which provides a more direct physical argument for a similar problem. I won't attempt to say too much about the physical interpretation of this statement; the interested reader is encouraged to check those two sources \cite{YujiRecipe,WWW}.

The proof I will present begins with a slight generalization of the one described in \cite[\S 2.7]{YujiRecipe}. Motivated by the physical interpretation, the result in \cite{YujiRecipe} is stated in terms of $n$-cocycles of $G$ with $n \geq 2$ and coefficients in $U(1) \subset \mathbb{C}^{\times}$, taken as a trivial $G$-module. The argument presented, however, works just as well for cocycles with coefficients in more general $G$-modules, so long as their image generates a finite subgroup. The technique is quite elementary, with roots in \cite{LHS}, the paper in which Hochschild and Serre introduced introduced the spectral sequence we'll encounter in the proof. For those more comfortable with the machinery of spectral sequences, there is a more high-tech proof of this statement given in Lemma 4.2.3 of \cite{Thibault} in the case where the coefficient module is torsion as an abelian group. I go a bit farther in this note by showing that any $n$-cocycle of $G$ with coefficients in \textit{any} $G$-module can be trivialized, provided $n \geq 3$.

In the next section, I'll recall some useful facts about the Lyndon-Hochschild-Serre (LHS) spectral sequence associated to an extension of a finite group $G$ by an abelian kernel. These facts will be used in the following section, in which I prove the stated result. In the Appendix, I include some background on group cohomology that I believe would be helpful to those less familiar with these topics. In particular, I give a brief introduction to derived functors, which readers may safely skip if they are not interested in the more abstract formalism. Then I recall some elementary facts and definitions about group cohomology itself.

\section{Facts about the LHS Spectral Sequence}
\label{sec:SpecSeq}

The main setting I am interested in is that of a group extension by an abelian kernel. Consider such an extension of $G$ by an abelian $A$ given by the short exact sequece $1 \rightarrow A \hookrightarrow \Gamma \twoheadrightarrow G \rightarrow 1$. It's a well-known fact (e.g. \cite[\S IV]{GroupCohomology}) that the isomorphism class of $\Gamma$ is fully determined by the data of an action of $G$ on $A$ and a class in $H^2(G;A)$ that I'll call the \textit{extension class}. Up to group isomorphism, $\Gamma$ is given by the set $A \times G$ and, if $c \in Z^2(G;A)$ represents the extension class, the group law is given by $(a,g) (b,h) = (a + g\cdot b + c(g,h), gh)$. I will denote this particular group by $A \rtimes^c G$ and when dealing with an extension $\Gamma \cong A \rtimes^c G$, I'll abuse notation by writing $\gamma = (a,g)$ for $\gamma \in \Gamma$.

A useful tool to analyze the cohomology of $\Gamma$ in terms of $G$ and $A$ is the Lyndon-Hochschild-Serre (LHS) spectral sequence. Introduced in \cite{LHS} for computing the cohomology of a quotient group more generally, it's quite nice in this specific setting. If $M$ is some $G$-module, then since $G \cong \Gamma / A$, it may also be considered as a $\Gamma$-module on which the normal subgroup $A$ acts trivially. Taking coefficients in $M$, the $E_2$-page LHS spectral sequence is of the form
\[
E_2^{r,s} = H^r(G;H^s(A;M))
\]
\noindent and the sequence converges to the cohomology of $\Gamma$ in the sense that $H^n(\Gamma;M) \cong \bigoplus_{r+s = n} E_{\infty}^{r,s}$. In this context, the action $G$\Large{$\curvearrowright$}\normalsize{$H^s(A;M)$} is induced by the actions $G$\Large{$\curvearrowright$}\normalsize{$A$} and $G$\Large{$\curvearrowright$}\normalsize{$M$}, and is defined at the level of cochains by 
\[
(g\cdot f)(a_1, \dots, a_s) = g \cdot f(g^{-1}\cdot a_1, \dots, g^{-1}\cdot a_s)
\]
Since $A$\Large{$\curvearrowright$}\normalsize{$M$} trivially, $H^1(A;M) = \operatorname{Hom_{Grp}}(A,M)$, and so the function $A \otimes H^1(A;M) \rightarrow M$ given by $a \otimes f \mapsto f(a)$ is a well-defined group homomorphism. Moreover, with the induced action on $H^1(A;M)$, this evaluation map is also a pairing of $G$-modules. As discussed in Appendix \ref{ap:C}, this then induces a cup product in cohomology of the following form.
\[
\cup_{\operatorname{ev}} \colon H^*(G;A) \otimes H^*(G;H^1(A;M)) \longrightarrow H^*(G;M)
\]

In order to prove the result claimed, I'll need to use a couple facts about the LHS spectral sequence. To orient ourselves, consider the following portion of the $E_2$-page corresponding to the short exact sequence $1 \rightarrow A \hookrightarrow \Gamma \twoheadrightarrow G \rightarrow 1$.
\begin{equation*}
\def\dx{2}
\def\dy{1}
\def\h{1}
\def\w{5}
\def\hbreak{2}
\def\vbreak{10}
\def\hoffset{-2}
\def\voffset{-1}
\def\hcenter{$n$}
\def\vcenter{$n$}
\def\vertax{$s$}
\def\horax{$r$}
\def\vertdots{1}
\def\hordots{1}
\begin{tikzpicture}[anchor=base, scale=0.8,every node/.style={scale=0.8}]

\path
(0,0) node[anchor=center] {$M^G$} ++(\dx,0) node[anchor=center] {$H^1(G)$} ++(\dx,0) node[anchor=center] {$\cdots$} ++(\dx,0) node[anchor=center] {$*$} ++(\dx,0) node[anchor=center] {$*$} ++(\dx,0) node[anchor=center] (target) {$H^n(G)$} ++(0.75*\dx,0) node[anchor=center] {$\cdots$}

(0,\dy)  node[anchor=center] {$H^1({A})^G$} ++(\dx,0) node[anchor=center] {$H^1(G;H^1({A}))$} ++(\dx,0) node[anchor=center] {$\cdots$} ++(\dx,0) node[anchor=center] (source) {$H^{n-2}(G;H^1({A}))$} ++(\dx,0) node[anchor=center] {$*$} ++(\dx,0) node[anchor=center] {$*$} ++ (0.75*\dx,0) node[anchor=center] {$\cdots$}

(0,1*\dy + 0.75*\dy) node[anchor=center] {\raisebox{6pt}{$\vdots$}} ++(\dx,0) node[anchor=center] {\raisebox{6pt}{$\vdots$}} ++(\dx,0) node[anchor=center] {\raisebox{6pt}{$\vdots$}} ++(\dx,0) node[anchor=center] {\raisebox{6pt}{$\vdots$}} ++(\dx,0) node[anchor=center] {\raisebox{6pt}{$\vdots$}} ++(\dx,0) node[anchor=center] {\raisebox{6pt}{$\vdots$}} ++(0.75*\dx,0) node[anchor=center] {\raisebox{6pt}{$\iddots$}} 
;

\ifnum \vertdots = 1
    \draw[->] (-0.5 * \dx, -\dy) -- (-0.5*\dx, \h * \dy + 1.25*\dy) node[anchor=south east] () {\vertax};
    \node[anchor=east] at (-0.5 * \dx,\h*\dy + 0.75*\dy) {\raisebox{6pt}{$\vdots$}};
    
\else
    \draw[->] (-0.5 * \dx, -\dy) -- (-0.5*\dx, \h * \dy + 0.5*\dy) node[anchor=south east] () {\vertax};
\fi
\ifnum \hordots = 1
    \draw[->] (-3/4 * \dx, -0.5*\dy) -- (\w * \dx + 1.25 * \dx, -0.5*\dy) node[anchor=north west] () {\horax};
    \node[anchor=north] at (\w*\dx + 0.75*\dx, -0.5 * \dy) {$\stackrel{}{\dots}$};
\else
    \draw[->] (-3/4 * \dx, -0.5*\dy) -- (\w * \dx + 0.5 * \dx, -0.5*\dy) node[anchor=north west] () {\horax};
\fi

\newcounter{dummy}
\setcounter{dummy}{\hoffset}

\foreach \i in {0,...,\w}{
    \ifnum \i < \hbreak
        \node[anchor=north] at (\i*\dx, -0.5 * \dy) {\i};
    \else    
        \ifnum \i = \hbreak
            \node[anchor=north] at (\i*\dx, -0.5 * \dy) {$\stackrel{}{\dots}$};
        \else
            \ifnum \inteval{\value{dummy}} < \inteval{0}
                \node[anchor=north] at (\i*\dx, -0.5 * \dy) {\hcenter\ \raisebox{1pt}{\scriptsize{$-$}} \normalsize{\fpeval{abs(\value{dummy})}}};
            \else
                \ifnum \inteval{\value{dummy}} = \inteval{0}
                    \node[anchor=north] at (\i*\dx, -0.5 * \dy) {\hcenter};
                \else
                    \node[anchor=north] at (\i*\dx, -0.5 * \dy) {\hcenter\ \raisebox{1pt}{\scriptsize{$+$}} \normalsize{\inteval{\value{dummy}}}};
                \fi
            \fi
            \stepcounter{dummy}
        \fi
    \fi
}

\setcounter{dummy}{\voffset}
\foreach \i in {0,...,\h}{
    \ifnum \i < \vbreak
        \node[anchor=east] at (-0.5 * \dx,\i*\dy) {\i};;
    \else    
        \ifnum \i = \vbreak
           \node[anchor=east] at (-0.5 * \dx,\i*\dy) {\raisebox{6pt}{$\vdots$}};
        \else
            \ifnum \inteval{\value{dummy}} < \inteval{0}
                \node[anchor=east] at (-0.5 * \dx,\i*\dy) {\vcenter\ \raisebox{1pt}{\scriptsize{$-$}} \normalsize{\fpeval{abs(\value{dummy})}}};
            \else
                \ifnum \inteval{\value{dummy}} = \inteval{0}
                    \node[anchor=east] at (-0.5 * \dx,\i*\dy) {\vcenter};
                \else
                    \node[anchor=east] at (-0.5 * \dx,\i*\dy) {\vcenter\ \raisebox{1pt}{\scriptsize{$+$}} \normalsize{\inteval{\value{dummy}}}};
                \fi
            \fi
            \stepcounter{dummy}
        \fi
    \fi
}

\draw[->] (source) -- (target) node[anchor=south,pos=0.5, rotate=-15] {$d_2$};
\draw[->] (target) -- ++(1.5*\dx,-0.75*\dy) node[anchor=south, pos=0.75, rotate=-15] {$d_2 = 0$};

\end{tikzpicture}
\end{equation*}

\noindent where $H^{\bullet}(-)$ is meant to be understood as $H^{\bullet}(-;M)$, and $(-)^G$ as $H^0(G;-)$, the $d_2$ arrows are the appropriate differentials, and $*$ indicates some other cohomolgy groups that aren't relevant. This is a first-quadrant spectral sequence since we take $E_2^{r,s} = 0$ when either $r<0$ or $s<0$. In general, differentials on the $E_k$-page are homomorphisms of the following form.
\[
d_k \colon E_k^{r,s} \longrightarrow E_k^{r+k,s+k-1}
\]
Hence in particular, the differentials leaving $E_k^{r,0}$ are all the zero map for any $k \geq 2$. By definition of page-turning, the bottom row of the sequence in subsequent pages is then given by the following quotient.
\begin{align*}
    E_{k+1}^{r,0} &= \frac{\operatorname{ker}\left(E_k^{r,0}\stackrel{d_k}{\rightarrow}0\right)}{\operatorname{im}\left(E_k^{r-k,k-1}\stackrel{d_k}{\rightarrow} E_k^{r,0}\right)}\\
    &= \frac{E_k^{r,0}}{\operatorname{im}\left(E_k^{r-k,k-1}\stackrel{d_k}{\rightarrow} E_k^{r,0}\right)}
\end{align*}
\noindent These quotients give us surjective homomorphisms $E_k^{r,0} \twoheadrightarrow E_{k+1}^{r,0}$ for all $k \geq 2$. Portions of the spectral sequence stabilize after the differentials with sources and targets within this portion are all the zero maps. By inspecting when the sources of the differentials hitting the bottom row are off the page, it's clear that $E_{r-1}^{r,0} = E_{\infty}^{r,0}$. Moreover, convergence of the LHS spectral sequence says that there is a homomorphism $E_{\infty}^{r,0} \hookrightarrow H^r(\Gamma; M)$. Altogether, there is the following composition.
\begin{equation}
\label{diag:LiftingComposition}
\begin{tikzcd}
	{H^r(G;M)=E_2^{r,0}} & \cdots & {E_{r-1}^{r,0}=E_{\infty}^{r,0}} & {H^r(\Gamma;M)}
	\arrow[two heads, from=1-1, to=1-2]
	\arrow[two heads, from=1-2, to=1-3]
	\arrow[hook, from=1-3, to=1-4]
\end{tikzcd}
\end{equation}
\noindent In \cite[\S III.1]{LHS}, the authors argue that this composition is given by the ``lifting homomorphism" $H^r(G;M) \rightarrow H^r(\Gamma;M)$ which is the map induced in cohomology by the surjective homomorphism $\Gamma \twoheadrightarrow G$. In particular, a class $\omega \in H^r(G;M)$ is lifted to a trivial class if it's in the kernel of one of these surjections, i.e. $\omega$ is in the image of a differential hitting the bottom row of the sequence.

The last thing I'll need from the LHS is a formula for the differentials on the $E_2$-page hitting the bottom row. These are homomorphisms of the form $d_2 \colon H^r(G;H^1(A;M)) \rightarrow H^{r+2}(G;M)$. Theorem 4 of \cite[\S III]{LHS} shows that, for any $\alpha \in H^r(G;H^1(A;M))$, its differential is given by $d_2(\alpha) = - \alpha \cup_{\operatorname{ev}}c$, where $\cup_{\operatorname{ev}}$ is the cup product induced by the evaluation pairing, and $c \in H^2(G;A)$ is the extension class. Using Definition \ref{def:CupProduct} (in Appendix \ref{ap:C}), we get the following equation  (in $M$) at the level of the representative cocycles.
\begin{equation}
\label{eq:DifferentialEquation}
    d_2(\alpha)(g_1,\dots,g_{r+2}) = - \alpha_{(g_1,\dots,g_r)}(g_1 \cdots g_r \cdot c(g_{r+1},g_{r+2}))
\end{equation}
\noindent My notation here is that $\alpha_{(g_1,\dots,g_{r})} \in H^1(A;M)$ is the image of the map $\alpha \colon G^{\times r} \rightarrow H^1(A;M)$.

\section{Trivializations}

Throughout this section, $G$ is some fixed finite group, $M$ a $G$-module, and $\omega \colon G^{\times n} \rightarrow M$ is a fixed $n$-cocycle of $G$ in $M$. In general, for a group homomorphism $\phi \colon K \rightarrow G$, I'll write $\phi^* \coloneqq H^i(\phi;M) \colon H^i(G;M) \rightarrow H^i(K;M)$ so that e.g. $\phi^*\omega = \omega \circ \phi^{\times n}$, while for a $G$-module homomorphism $f\colon M \rightarrow N$, I'll write $f_*\coloneqq H^i(G;M) \rightarrow H^i(G;N)$ so that e.g. $f_* \omega = f \circ \omega$. The groups or modules may change, but the pattern will remain the same.

\subsection{Main results}

The image of $\omega$ generates a subgroup $\Omega \coloneqq \langle \operatorname{im}(\omega)\rangle \subseteq M$. The claim I will prove is that, if $\Omega$ is finite and $n \geq 2$, then $\omega$ may be trivialized in a finite extension of $G$ by an abelian kernel. More precisely, there is the following theorem.

\begin{Theorem}
\label{thm:Main}
If $|\Omega| < \infty$, then there exists a short exact sequence of finite groups
    \begin{equation}
    \label{diag:FiniteExt}
    \begin{tikzcd}
        1 & {A} & {\Gamma} & G & 1
        \arrow[from=1-1, to=1-2]
        \arrow["{\iota}", hook, from=1-2, to=1-3]
        \arrow["{\pi}", two heads, from=1-3, to=1-4]
        \arrow[from=1-4, to=1-5]
    \end{tikzcd}
    \end{equation}
    and a cochain $\alpha \in C^{n-1}(G;M)$ such that $A$ is abelian and $\pi^* \omega = \delta_{\Gamma} \alpha \in Z^n(\Gamma;M)$. In particular, $\pi^*\omega = 0 \in H^n(\Gamma;M)$
\end{Theorem}

\begin{Remark}
Given the extension (\ref{diag:FiniteExt}), the inclusion $\iota \colon A \hookrightarrow \Gamma$ induces a cochain map $\iota^* \colon C^{\bullet}(\Gamma;M) \rightarrow C^{\bullet}(A;M)$. If $\pi^*\omega$ is trivialized by $\alpha$, then  $\iota^*\alpha$ is necessarily a cocycle of $A$. Indeed, it's not hard to check
\begin{align*}
    \delta_A (\iota^* \alpha) &= \iota^* (\delta_{\Gamma} \alpha)\\
    &= \iota^*(\pi^* \omega)\\
    &= (\pi \circ \iota)^*\omega\\
    &= 0
\end{align*}
\end{Remark}

\noindent Notice that since $|G^{\times n}| = |G|^n < \infty$ and $M$ is abelian, $\Omega$ is a finitely-generated abelian group. Thus, the fundamental theorem of finitely generated abelian groups implies that $|\Omega| < \infty$ if and only if $\Omega$ is torsion. If $M$ is itself torsion, then $\Omega$ is guaranteed to be finite and, in this way, the above theorem immediately implies the following.

\begin{Cor}
\label{cor:Main}
If $G$ is a finite group and $M$ a $G$-module that is torsion as an abelian group, then any $n$-cocycle of $G$ in $M$ for $n \geq 2$ is trivializable in an extension of $G$ by a finite abelian group.
\end{Cor}

A more direct proof of Corollary \ref{cor:Main} that bypasses Theorem \ref{thm:Main} is given as Lemma 4.2.3 of \cite{Thibault} using some general spectral sequence machinery applied to the relevant LHS spectral sequence. I'll prove Theorem \ref{thm:Main} in the next section. In the meantime, assuming it's true, we also get the following theorem.

\begin{Theorem}
Let $M$ be any $G$-module, and $\omega$ an $n$-cocycle of $G$ in $M$ with $n \geq 3$, then there is a finite extension $\widetilde{\Gamma}$ of $G$ in which $\omega$ is trivializable.
\end{Theorem}
\begin{Proof}
Define $M_T$ to be the $\Z$-torsion submodule of $M$, i.e. the sub $G$-module given by all elements of $M$ of finite order. Now consider the short exact sequence of trivial $G$-modules (i.e. abelian groups) given by $0 \rightarrow \Z \hookrightarrow \Q \twoheadrightarrow \bigslant{\Q}{\Z} \rightarrow 0$. Writing $\otimes = \otimes_{\Z}$, we can apply the right exact functor $M \otimes -$ to this sequence and compute the kernel of the left-most map. It's not hard to show that this is precisely $M_T$, and so we have the four-term exact sequence $0 \rightarrow M_T \hookrightarrow M \rightarrow M \otimes \Q \twoheadrightarrow M \otimes \bigslant{\Q}{\Z} \rightarrow 0$. This decomposes into two short exact sequences as follows.
\[
\begin{tikzcd}[sep=small]
	&& 0 && 0 \\
	&&& {M/M_T} \\
	0 & {M_T} & M && {M\otimes \mathbb{Q}} & {M\otimes \mathbb{Q}/\mathbb{Z}} & 0
	\arrow[from=1-3, to=2-4]
	\arrow[from=2-4, to=1-5]
	\arrow[from=2-4, to=3-5]
	\arrow[from=3-1, to=3-2]
	\arrow[from=3-2, to=3-3]
	\arrow[from=3-3, to=2-4]
	\arrow[from=3-3, to=3-5]
	\arrow[from=3-5, to=3-6]
	\arrow[from=3-6, to=3-7]
\end{tikzcd}
\]

Next, since the (normalized) bar resolution is given by $i$-dimensional \textit{free} $\Z[G]$-modules $\overline{\Z G}_i$ (see Appendix \ref{sec:GModules}), the covariant functors $\operatorname{Hom}_{G}(\overline{\Z G}_i,-)$ are all exact. We therefore have two corresponding short exact sequences of cochain complexes. The first is $0 \rightarrow C^{\bullet}(G;\bigslant{M}{M_T}) \rightarrow C^{\bullet}(G;M\otimes \Q) \rightarrow C^{\bullet}(G;M\otimes \bigslant{\Q}{\Z}) \rightarrow 0$ and yields the following long exact sequence in cohomology.
\begin{equation}
\label{eq:LES1}
\begin{tikzpicture}
    \node (c1) at (0,0) {$H^i\left(G;\bigslant{M}{M_T}\right)$};
    \node[anchor=west] (c2) at ($(c1.east) + (1,0)$) {$H^i(G;M\otimes Q)$};
    \node[anchor=west] (c3) at ($(c2.east) + (1,0)$) {$H^i(G;M \otimes \bigslant{\Q}{\Z})$};
    \node[anchor=east] (t1) at ($(c3.east) + (0,1)$) {$H^{i-1}\left(G;M\otimes \bigslant{\Q}{\Z}\right)$};
    \node[anchor=east] (tdots) at ($(t1.west) + (-0.5,0)$) {$\cdots$};
    \node[anchor=west] (b1) at ($(c1.west) + (0,-1)$) {$H^{i+1}\left(G;\bigslant{M}{M_T}\right)$};
    \node[anchor=west] (bdots) at ($(b1.east) + (0.5,0)$) {$\cdots$};
    \draw[->] (tdots.east) -- (t1.west);
    \draw[->] (b1.east) -- (bdots.west);
    \draw[->] (c1.east) -- (c2.west);
    \draw[->] (c2.east) -- (c3.west);
    \draw[->] (c3.east) --++ (0.5,0) --++ (0,-0.5) -- ($(c1.west) + (-0.5,-0.5)$) --++ (0,-0.5) -- (b1.west);
    \draw[->] (t1.east) --++ (0.5,0) --++ (0,-0.5) -- ($(c1.west) + (-0.5,0.5)$) --++ (0,-0.5) -- (c1.west);
\end{tikzpicture}
\end{equation}
Proposition 9.5 in \cite[\S III.9]{GroupCohomology} says that for any subgroup $K \subseteq G$ and $G$-module $N$, there are homomorphisms $H^{\bullet}(G;N) \rightarrow H^{\bullet}(G;N)$ given by $c \mapsto [G:K]c$ and that this factors through $H^{\bullet}(K;N)$. Taking $K = \{1\}$ thus implies that $|G|c = 0$ for any class $c \in H^{\bullet}(G;N)$. If $N$ is a vector space over some field containing $|G|$, then $H^{\bullet}(G;N) = \{0\}$ because $c = (|G|^{-1}|G|)c=|G|^{-1}(|G|c)=|G|^{-1}(0) = 0$ for any $c$. It's not hard to show that $M \otimes \Q$ is indeed a vector space over $\Q$, and thus that $H^i(G;M\otimes \Q) = 0$ for all $i$. The exactness of (\ref{eq:LES1}) then implies that $H^i(G;\bigslant{M}{M_T}) \cong H^{i-1}(G;M\otimes \bigslant{\Q}{\Z})$.

Let $p \colon M \rightarrow \bigslant{M}{M_T}$ be the quotient map, then $p\circ\omega = p_*\omega \in H^n(G; \bigslant{M}{M_T})$. Since $M \otimes \bigslant{\Q}{\Z}$ is torsion as an abelian group, Corollary \ref{cor:Main} implies that if $n-1 \geq 2$, then there is a finite extension $\pi \colon \Gamma \twoheadrightarrow G$ of $G$ such that $\pi^*(p_*\omega) = 0 \in H^{n}(\Gamma; \bigslant{M}{M_T}) \cong H^{n-1}(\Gamma; M \otimes \bigslant{\Q}{\Z})$. The other short exact sequence of chain complexes we'll need to consider is $0 \rightarrow C^{\bullet}(\Gamma; M_T) \rightarrow C^{\bullet}(\Gamma; M) \rightarrow C^{\bullet}(\Gamma; M/M_T) \rightarrow 0$. A snapshot of the associated long exact sequence in cohomology looks like the following.
\begin{equation}
\label{eq:LES2}
\begin{tikzpicture}
    \node (c1) at (0,0) {$H^i\left(\Gamma;M_T\right)$};
    \node[anchor=west] (c2) at ($(c1.east) + (0.75,0)$) {$H^i(\Gamma;M)$};
    \node[anchor=west] (c3) at ($(c2.east) + (0.75,0)$) {$H^i(\Gamma;M/M_T)$};
    \draw[->] (c1.east) -- (c2.west) node[pos=0.5,anchor=south] {$j_*$};
    \draw[->] (c2.east) -- (c3.west) node[pos=0.5,anchor=south] {$p_*$};
    \draw[->] (c3.east) --++ (0.25,0) node[pos=1, anchor=west] {$\cdots$};
    \draw[<-] (c1.west) --++ (-0.25,0) node[pos=1, anchor=east] {$\cdots$};
\end{tikzpicture}
\end{equation}
By assumption, $\pi^*(p_*\omega) = 0$, but we also have $\pi^*(p_*\omega) = p \circ \omega \circ \pi^{\times n} = p_*(\pi^*\omega)$. This means that $\pi^*\omega \in \operatorname{ker}(p_*)$ and thus, by the exactness of (\ref{eq:LES2}), there is a cocycle $\beta \colon \Gamma^n \rightarrow M_T$ such that $j_* \beta = \pi^*\omega$, with $j \colon M_T \hookrightarrow M$ the inclusion.

In particular we have $n \geq 2$ and moreover, $j_*\beta$ is an $n$-cocycle of $\Gamma$ in $M$ whose image generates a finite group. This is precisely the hypothesis of Theorem \ref{thm:Main}, and so there will be a finite extension $\widetilde{\pi}\colon\widetilde{\Gamma} \twoheadrightarrow \Gamma$ of $\Gamma$ such that $\widetilde{\pi}^*(j_*\beta) = 0 \in H^n(\widetilde{\Gamma};M)$. In turn, this tells us that our original cocycle is also trivialized here too, since $\widetilde{\pi}^*(j_*\beta) = \widetilde{\pi}^*(\pi^*\omega) = (\pi \circ \widetilde{\pi})^*\omega$.
\end{Proof}

\subsection{The Proof of Theorem \ref{thm:Main}}
\label{sec:Proof}

Following \cite[\S 2.7]{YujiRecipe}, the strategy will be to first trivialize $\omega$ in an extension of $G$ by a free abelian (hence infinite) kernel $\mathcal{A}$. Next, I then tensor $\mathcal{A}$ with a finite cyclic group to get $A$. The reason we need to assume that $|\Omega| < \infty$ is in order to actually construct the finite cyclic group in question.

Define $\mathcal{A}$ to be the free abelian group generated by symbols $\{\times_{g,h}\ |\ g, h \neq 1 \in G\}$. Alternatively, this may be viewed as first being freely generated by $G \times G$, and then quotiented by a relation setting $\times_{1,g} = \times_{g,1} = 0$ for any $g \in G$. The reason for this is that $\mathcal{A}$ should be the universal target for normalized 2-cocycles of $G$. In order to accomplish this, the map $c \colon G\times G \rightarrow \mathcal{A}$ given by $c(g,h)=\times_{g,h}$ should itself be a 2-cocycle. For this to make sense, we equip $\mathcal{A}$ with a $G$-action defined by
\[
    g \cdot \times_{h,k} = \times_{gh,k} - \times_{g,hk} + \times_{g,h}
\]
\noindent for any $g,h,k \in G$. This is ``just-so" so that $c$ satisfies the 2-cocycle condition (i.e. $\delta c = 0$). This $G$-module $\mathcal{A}$ then represents the $\operatorname{Ab}$-enriched functor of taking 2-cocycles of $G$.
\[
Z^2(G;-) \colon G\operatorname{Mod} \longrightarrow \operatorname{Ab}
\]
Its universal property is then as follows: for any $G$-module $M$ and any normalized $2$-cocycle $\mu \colon G \times G \rightarrow M$ there is a unique $G$-module homomorphism $\varphi \colon \mathcal{A} \rightarrow M$ making the following the following diagram commute

\begin{equation*}
    \label{diag:UniversalCocycle}
    \begin{tikzcd}
    	{G \times G} & {\mathcal{A}} \\
    	& M
    	\arrow["{{c}}", from=1-1, to=1-2]
    	\arrow["\mu"', from=1-1, to=2-2]
    	\arrow["\varphi", dashed, from=1-2, to=2-2]
    \end{tikzcd}
\end{equation*}

The $G$-module $\mathcal{A}$ and the class in $H^2(G;\mathcal{A})$ represented by $c$ determines an extension $\mathcal{G} \cong \mathcal{A} \rtimes^c G$ of $G$ by $\mathcal{A}$. This is described by the following short exact sequence.
\[\begin{tikzcd}
    1 & \mathcal{A} & {\mathcal{G}} & G & 1
    \arrow[from=1-1, to=1-2]
    \arrow["{j}", hook, from=1-2, to=1-3]
    \arrow["{p}", two heads, from=1-3, to=1-4]
    \arrow[from=1-4, to=1-5]
\end{tikzcd}\]
The next step is to show that $\omega$ is trivial as an $n$-cocycle of ${\mathcal{G}}$ in $M$, which means that $p^*\omega$ represents $0 \in H^n({\mathcal{G}};M)$. Here, $p^*$ is the lifting homomorphism, so it's given by the composition (\ref{diag:LiftingComposition}). Since $\omega$ represents a class in $H^n(G;M) = E_2^{n,0}$, it would suffice to show that, as the pages of the spectral sequence are turned, $\omega$ eventually gets sent to the kernel of one of the surjections in the composition. Fortunately, $\omega$ is already killed by the first surjection $E_{2}^{n,0} \twoheadrightarrow E_3^{n,0}$. By definition, this is the same as being in the image of the $d_2$ differential, i.e. that there exists some $b \in Z^{n-2}(G;H^1(\mathcal{A};M))$ representing a class in $H^{n-2}(G;H^1(\mathcal{A};M)) = E_2^{n-2,1}$ such that
\begin{align*}
\begin{split}
    \omega(g_1,\dots,g_n) &= \left(d_2 b\right)(g_1,\dots,g_n)\\
    &= -b_{(g_1,\dots,g_{n-2})} \left(g_1 \cdots g_{n-2} \cdot {c}(g_{n-1},g_n)\right)\\
\end{split}
\end{align*}
\noindent where I used the formula (\ref{eq:DifferentialEquation}) for the differential $d_2$. Incidentally, this is where the requirement that $n \geq 2$ comes into play; if $n = 0$ or $n = 1$, then there are no differentials hitting $E_k^{n,0}$ for \textit{any} $k \geq 2$. I'll show that such a $b$ exists by directly constructing it, and I'll do so in two steps. First, let $Y$ be a $(n-2)$-cochain of $\mathcal{G}$ with values in $H^1(\mathcal{A};M) = \operatorname{Hom_{Grp}}(\mathcal{A},M)$ that is defined as follows. 
\begin{equation}
\label{eq:YCochain}
    Y_{(g_1,\dots,g_{n-2})}(\times_{g_{n-1},g_n}) = - \omega(g_1,\dots,g_n)
\end{equation} 
\noindent This gives a well-defined cochain with values in $H^1(\mathcal{A};M)$ by $\Z$-linearly extending to all of $\mathcal{A}$. Then I define a new cochain by slightly modifying $Y$ in the following manner.
\[
b_{(g_1, \dots, g_{n-2})} \coloneqq \left((g_1 \cdots g_{n-2})^{-1}\right) \cdot Y_{(g_1,\dots g_{n-2})}
\]
Therefore, by construction, $d_2(b) = \omega$. All that remains to check is that $b$ is in fact a cocycle. I'll show this by using the Leibniz rule (\ref{eq:Leibniz}), and leveraging the fact that both $\omega$ and $c$ are cocycles of $G$.
\begin{align*}
0 &= -\delta_G\omega(g_1,\dots,g_{n+1})\\
&= \delta_G(b\cup_{\operatorname{ev}}{c})(g_1,\dots,g_{n+1})\\
&=\big(\delta_G b \cup_{\operatorname{ev}} {c} + (-1)^{n-2} b\cup_{\operatorname{ev}}\delta_G{c}\big)(g_1,\dots,g_{n+1})\\
&=\big(\delta_G b\cup_{\operatorname{ev}} c \big)(g_1,\dots,g_{n+1})\\
&= (\delta_G b)_{(g_1, \dots, g_{n-1})}\big(g_1 \cdots g_{n-1} \cdot c(g_n,g_{n+1})\big)\\
&= (g_1 \cdots g_{n-1})  \cdot \big((g_1 \cdots g_{n-1})^{-1} \cdot (\delta_G b)_{(g_1, \cdots, g_{n-1})}\big)(c(g_n, g_{n+1}))\\
&= \big((g_1 \cdots g_{n-1})^{-1} \cdot (\delta_G b)_{(g_1, \cdots, g_{n-1})}\big)(c(g_n, g_{n+1}))
\end{align*}
where in the final step, I acted on both sides by $(g_1 \cdots g_{n-1})^{-1}$, via the $G$-action on $M$. Now, since $c(g_n, g_{n+1}) = \times_{g_n,g_{n+1}}$ are \textit{generators} of $\mathcal{A}$, the above can be viewed as an equation in $H^1(\mathcal{A};M) = \operatorname{Hom_{Grp}}(\mathcal{A}, M)$.
\[
0 = (g_1 \cdots g_{n-1})^{-1} \cdot (\delta_G b)_{(g_1, \cdots, g_{n-1})}
\]
I again act on both sides, this time by $g_1 \cdots g_{n-1}$ via the $G$-action on $H^1(\mathcal{A};M)$, which finally gives us $\delta_G b = 0 \in C^{n+1}(G;H^1(\mathcal{A};M)$. Thus $b$ is indeed a cocycle as required, and we have that $p^*\omega = 0 \in H^n(\mathcal{G};M)$.

At this point, our task is almost complete. The only thing left to do now is to make everything finite. This can be done via the exact same argument and using the same formulas, but replacing $\mathcal{A}$ with the finite group $A \coloneqq \mathcal{A}\otimes_{\Z}\Z/N$ for some appropriately large natural number $N$. The reason that $N$ needs to be sufficiently large is so that the formula (\ref{eq:YCochain}) defining the cochain $Y$, and hence the cocycle $b$ remains well-defined. In particular, $Y$ and $b$ need to actually take values in $H^1(A;M) = \operatorname{Hom_{Grp}}(A, M)$. Consider a case where $N$ is smaller than the order of $\omega(h_1,\dots,h_n)$ for some choice of $h_i$'s in $G$, and denote the generator of $\Z/N$ by $1_N$. In this case, there is the following contradiction
\begin{align*}
N Y(h_1,\dots,h_{d-2})(\times_{h_{d-1},h_d} \otimes 1_N) &= -N \omega(h_1,\dots,h_d) \neq 0\\
&=Y(h_1,\dots,h_{d-2})(N \times_{h_{d-1},h_d} \otimes 1_N) = 0
\end{align*}
\noindent Since I assume that $\Omega$ is finite (and hence torsion), one way to guarantee that this doesn't happen is to take $N$ to be the least common multiple of the orders of non-zero elements in $\operatorname{im}(\omega)$. In this way, the above contradiction is avoided since, for any collection of positive integers $\{n_i\}_{i=1}^k$, $\operatorname{lcm}(n_1,\dots,n_k) \geq n_i$ for any $i$. Moreover, $N\omega(g_1,\dots,g_n) = 0$ for any choice of $g_i$'s in $G$, as is required. Replacing $\mathcal{A} \cong \Z^{(|G|-1)^2}$ with $A \cong (\Z/N)^{(|G|-1)^2}$, the exact same argument then proves the theorem.

\hfill$\square$

\appendix

\section*{Appendix}

There's no shortage of things to say about group cohomology. Indeed, it's a concept that leads itself to interpretations in a variety of settings. As a consequence, there are many equivalent ways to define the cohomology of some group $G$ with coefficients in some module $M$. The most direct definition is in terms of a derived functor, so I'll first spend some time saying what those are, as well as a bit about \textit{why} they are what they are.

\section{Chain homotopy theory}
\label{sec:ChainHomotopy}
As I'll say in the next section, the cohomology of a group $G$ with coefficients in some $G$-module is given by a certain derived functor. When I initially learned about cohomology, and derived functors more generally, they seemed quite arcane and computing them felt like following a random sequence of steps. The purpose of this section is to slightly demystify this concept by highlighting their homotopical nature. In modern language a ``homotopy theory" refers to (equivalence classes of) $(\infty,1)$-categories. These can be presented 1-categorically via \textit{model categories}, which is the context I'll use. 

A model category consists of a complete and cocomplete 1-category $\mathcal{M}$ together with three distinguished classes of morphisms called weak equivalences, fibrations, and cofibrations. These must then satisfy a bunch of axioms I won't spell out here (see e.g. \cite{Riehl, DwyerSpalinski, Hovey}), but that seek to abstract the homotopy theory studied in algebraic topology. In particular, it's common to consider spaces and continuous maps only ``up to homotopy", and this intuition is captured by the \textit{homotopy category} associated to a model category. This is obtained by formally inverting (i.e. localizing at) the weak equivalences. For a model category $\mathcal{M}$ with weak equivalences $\mathcal{W}$, I'll denote its canonical localization functor by $\operatorname{loc}_{\mathcal{M}} \colon \mathcal{M} \rightarrow \mathcal{M}\left[\mathcal{W}^{-1}\right] \coloneqq \text{Ho}\mathcal{M}$. I will also denote the fibrant and cofibrant replacement functors by $\mathcal{F}ib$ and $\mathcal{C}of$ respectively, and they come with the following natural isomorphisms.
\[
\begin{tikzcd}
	{\mathcal{M}} & {\text{Ho}\mathcal{M}} & {\mathcal{M}} & {\text{Ho}\mathcal{M}}
	\arrow[""{name=0, anchor=center, inner sep=0}, "{\mathcal{C}of}", curve={height=-12pt}, from=1-1, to=1-2]
	\arrow[""{name=1, anchor=center, inner sep=0}, "{\text{loc}_\mathcal{M}}"', curve={height=12pt}, from=1-1, to=1-2]
	\arrow[""{name=2, anchor=center, inner sep=0}, "{\mathcal{F}ib}", curve={height=-12pt}, from=1-3, to=1-4]
	\arrow[""{name=3, anchor=center, inner sep=0}, "{\text{loc}_\mathcal{M}}"', curve={height=12pt}, from=1-3, to=1-4]
	\arrow["\wr", shorten <=3pt, shorten >=3pt, Rightarrow, from=0, to=1]
	\arrow["\wr", shorten <=3pt, shorten >=3pt, Rightarrow, from=3, to=2]
\end{tikzcd}
\]

Let $F \colon \mathcal{M} \rightarrow \mathcal{K}$ be some (arbitrary) functor between model categories. $F$ may or may not preserve the weak equivalences, but if it does, then it also induces a functor between the corresponding homotopy categories. If $F$ doesn't preserve weak equivalences, then we can try to universally approximate it by one that does. This is the job of a derived functor.

\begin{Definition}
If they exist, the \textit{total right derived functor} $\mathscr{R}F$ of $F$ is a \textit{left} Kan extension while the \textit{total left derived functor} $\mathscr{L}F$ of $F$ is a \textit{right} Kan extension, both extending $\operatorname{loc}_{\mathcal{K}} \circ F$ along $\operatorname{loc}_{\mathcal{M}}$. This is summarized by the following diagrams.
\[
\begin{tikzcd}[column sep=tiny,row sep=scriptsize]
	{\mathcal{M}} && {\operatorname{Ho}\mathcal{K}} && {\mathcal{M}} && {\operatorname{Ho}\mathcal{K}} \\
	& {\text{Ho}\mathcal{M}} &&&& {\text{Ho}\mathcal{M}}
	\arrow[""{name=0, anchor=center, inner sep=0}, "\operatorname{loc}_{\mathcal{K}} \circ F", from=1-1, to=1-3]
	\arrow["{{\text{loc}_\mathcal{M}}}"', from=1-1, to=2-2]
	\arrow[""{name=1, anchor=center, inner sep=0}, "\operatorname{loc}_{\mathcal{K}} \circ F", from=1-5, to=1-7]
	\arrow["{{\text{loc}_\mathcal{M}}}"', from=1-5, to=2-6]
	\arrow["{{\mathscr{R}F}}"', dashed, from=2-2, to=1-3]
	\arrow["{{\mathscr{L}F}}"', dashed, from=2-6, to=1-7]
	\arrow[shorten <=3pt, Rightarrow, from=0, to=2-2]
	\arrow[shorten >=3pt, Rightarrow, from=2-6, to=1]
\end{tikzcd}
\]
\end{Definition}

\noindent There are certain functors, called left (resp. right) Quillen functors, that always admit left (resp. right) derived functors. In particular, A cocontinuous functor that preserves cofibrations as well as the subset of cofibrations that are weak equivalences is a left Quillen functor, while a continuous functor that preserves fibrations in the same manner is a right Quillen functor. By Corollary 4.2.4 of \cite{Riehl}, we may compute these derived functors as $\mathscr{L}F = F \circ \mathcal{C}of$ and $\mathscr{R}F = F \circ \mathcal{F}ib$. These are well-defined functions on objects since a (co)fibrant object is, in particular, an object of $\mathcal{M}$. On morphisms, these are only well-defined up to a weak equivalence, hence why they define functors to $\operatorname{Ho}\mathcal{K}$.

Now consider some abelian category $\mathcal{C}$. The category $\text{Ch}(\mathcal{C})$ of chain complexes internal to $\mathcal{C}$ has two full subcategories given by complexes bounded above and below respectively. Since we can identify any chain complex with a cochain complex by simply reflecting the indices, i.e. by $X_{\bullet} \mapsto X^{-\bullet}$, I'll write $\operatorname{Ch}_{\bullet}(\mathcal{C})$ for chain complexes bound below and $\operatorname{Ch}^{\bullet}(\mathcal{C})$ for chain complexes bound above, which are the same as cochain complexes bound below. These have some standard and well-studied model structures, given as follows.
\begin{Facts} Let $\mathcal{C}$ be an abelian category,  $\operatorname{Ch}_{\bullet}(\mathcal{C})$, and  $\operatorname{Ch}^{\bullet}(\mathcal{C})$ the (co)chain complexes in $\mathcal{C}$ bound from below.
\begin{enumerate}
    \item The \textit{projective model structure on $\mathcal{C}$}, defined on $\operatorname{Ch}_{\bullet}(\mathcal{C})$, is given by:
    \begin{itemize}
        \item The fibrations are epimorphisms $X_n \twoheadrightarrow Y_n$ in each degree.
        \begin{itemize}
            \item The unique chain map $X_{\bullet} \rightarrow 0_{\bullet}$ is a fibration, thus every object is fibrant.
        \end{itemize}
        \item The cofibrations are monomorphisms $X_n \hookrightarrow Y_n \twoheadrightarrow P_n$ with projective cokernels in each degree.
        \begin{itemize}
            \item The unique chain map $0_{\bullet} \rightarrow X_{\bullet}$ is a cofibration if and only if each $X_n$ is projective.
        \end{itemize}
        \item Right exact functors $F \colon \mathcal{C} \rightarrow \mathcal{D}$ induce left Quillen functors with respect to the projective model structure.
    \end{itemize}
        \item The \textit{injective model structure on $\mathcal{C}$}, defined on $\operatorname{Ch}^{\bullet}(\mathcal{C})$, is given by:

    \begin{itemize}
        \item The fibrations are epimorphisms $J^n\hookrightarrow X^n \twoheadrightarrow Y^n$ with injective kernels in each degree.
        
        \begin{itemize}
            \item The unique cochain map $X^{\bullet} \rightarrow 0^{\bullet}$ is a fibration if and only if each $X^n$ is injective.
        \end{itemize}
        
        \item The cofibrations are monomorphisms $X^n \hookrightarrow Y^n$ in each degree.
    
        \begin{itemize}
            \item The unique cochain map $0^{\bullet} \rightarrow X^{\bullet}$ is a cofibration, thus every object is cofibrant.
        \end{itemize}
        \item Left exact functors $F \colon \mathcal{C} \rightarrow \mathcal{D}$ induce right Quillen functors with respect to the injective model structure.
    \end{itemize}
    
    \item The weak equivalences in both model structures are \textit{quasi-isomorphisms}, $\mathcal{Q}$. These are (co)chain maps inducing isomorphisms in (co)homology in each degree. The \textit{derived categories of $\mathcal{C}$} are the homotopy categories.
    \[
    \begin{tikzcd}[row sep=tiny]
    	{\mathscr{D}_{\bullet}(\mathcal{C}) \coloneqq \operatorname{Ch}_{\bullet}(\mathcal{C}) \left[\mathcal{Q}^{-1}\right]\simeq \operatorname{Ch}_{\bullet}(\mathcal{P_C}) \left[\mathcal{Q}^{-1}\right]} \\
    	{\mathscr{D}^{\bullet}(\mathcal{C}) \coloneqq \operatorname{Ch}^{\bullet}(\mathcal{C}) \left[\mathcal{Q}^{-1}\right]\simeq \operatorname{Ch}^{\bullet}(\mathcal{J_C}) \left[\mathcal{Q}^{-1}\right]}
    \end{tikzcd}
    \]
    where $\mathcal{P_C}$ and $\mathcal{J_C}$ are the full subcategories of projective and injective objects of $\mathcal{C}$ respectively.
\end{enumerate}
\end{Facts}
The projective model structure is discussed in \cite[\S 7]{DwyerSpalinski} in the context of $\mathcal{C} = R \operatorname{Mod}$ for a ring $R$. The injective model structure for \textit{unbounded} chain complexes in $R\operatorname{Mod}$ is discussed in \cite[\S 2.3]{Hovey}. The \href{https://ncatlab.org/nlab/show/model+structure+on+chain+complexes}{nLab webpage} \cite{nLab} collects many of these facts together in one place, along with a long list of references for the interested reader.

In homological algebra, what one typically computes are what I'll refer to as \textit{$n^{th}$ derived functors}. For example, the $n^{th}$ right derived functor $\mathscr{R}^nF$ of a left exact $F \colon \mathcal{C} \rightarrow \mathcal{D}$ is the following composition.
\[
\begin{tikzcd}
	{\mathcal{C}} & {\text{Ch}^{\bullet}(\mathcal{C})} & {\text{Ch}^{\bullet}(\mathcal{J_C})\left[\mathcal{Q}^{-1}\right]} & {\mathscr{D}^{\bullet}(\mathcal{D})} & {\mathcal{D}}
	\arrow[hook, from=1-1, to=1-2]
	\arrow["{\mathcal{F}ib}", from=1-2, to=1-3]
	\arrow["F", from=1-3, to=1-4]
	\arrow["{{H^n}}", from=1-4, to=1-5]
\end{tikzcd}
\]
where the embedding $\mathcal{C} \hookrightarrow \text{Ch}^{\bullet}(\mathcal{C})$ is given by placing an object in degree zero. Fibrantly replacing this complex gives precisely an injective resolution of the object. The similar situation in the projective model structure leads to a projective resolution of the object. After evaluating the total derived functor, we then read off the its (co)homology objects. Anyway, I'll conclude with some well-known and important examples of derived functors.
\begin{Examples}
\label{ex:DerivedFuncs}Let $\mathcal{C}$ be a symmetric monoidal abelian category.
\begin{enumerate}
    \item $\otimes\colon \operatorname{Ch}(\mathcal{C} \times \mathcal{C}) \rightarrow \operatorname{Ch}(\operatorname{Ab})$ is a right exact functor, and thus admits a (total) left derived functor $\mathscr{L}\otimes$ with respect to the projective model structure on $\operatorname{Ch}_{\bullet}(\mathcal{C}\times \mathcal{C}) = \operatorname{Ch}_{\bullet}(\mathcal{C}) \times \operatorname{Ch}_{\bullet}(\mathcal{C})$. We then obtain the familiar $\mathscr{L}^n\otimes(X,Y) = \operatorname{Tor}_{\mathcal{C}}^n(X,Y)$. Moreover, using the data available from the definitions, we can build the following 2-cell.
\[
\begin{tikzcd}
	{\operatorname{Ch}_{\bullet}(\mathcal{C})\times \operatorname{Ch}_{\bullet}(\mathcal{C})} && {\operatorname{Ch}_{\bullet}(\mathcal{C})} \\
	{\mathscr{D}_{\bullet}(\mathcal{C}) \times \mathscr{D}_{\bullet}(\mathcal{C})} && {\mathscr{D}_{\bullet}(\mathcal{C})}
	\arrow[""{name=0, anchor=center, inner sep=0,pos=.47}, "\otimes", from=1-1, to=1-3]
	\arrow[""{name=1, anchor=center, inner sep=0}, "{\operatorname{loc}}", curve={height=-12pt}, from=1-1, to=2-1]
	\arrow[""{name=2, anchor=center, inner sep=0}, "{\mathcal{C}of_{\mathcal{C}}}"', curve={height=12pt}, from=1-1, to=2-1]
	\arrow[""{name=3, anchor=center, inner sep=0}, "{\operatorname{loc}}"', curve={height=12pt}, from=1-3, to=2-3]
	\arrow[""{name=4, anchor=center, inner sep=0}, "{\mathcal{C}of_{\mathcal{C}}}", curve={height=-12pt}, from=1-3, to=2-3]
	\arrow[""{name=5, anchor=center, inner sep=0}, "{\mathscr{L}\otimes}"', from=2-1, to=2-3]
	\arrow["\sim", shorten <=5pt, shorten >=5pt, Rightarrow, from=2, to=1]
	\arrow["\sim"', shorten <=5pt, shorten >=5pt, Rightarrow, from=4, to=3]
	\arrow["",shorten <=5pt, shorten >=5pt, Rightarrow, from=5, to=0]
\end{tikzcd}
\]
In particular, this shows that the cofibrant replacement functor $\mathcal{C}of_{\mathcal{C}}$ is lax monoidal.

\item $\operatorname{Hom}_{\mathcal{C}} \colon \operatorname{Ch}(\mathcal{C}^{\operatorname{op}} \times \mathcal{C}) \rightarrow \operatorname{Ch}(\operatorname{Ab})$ is left exact and so admits a right derived functor $\mathscr{R}\operatorname{Hom}_{\mathcal{C}}$ with respect to the injective model structure on $\operatorname{Ch}^{\bullet}(\mathcal{C}^{\operatorname{op}} \times \mathcal{C}) = \operatorname{Ch}_{\bullet}(\mathcal{C})^{\operatorname{op}} \times \operatorname{Ch}^{\bullet}(\mathcal{C})$. We then have $\mathscr{R}^n\operatorname{Hom}_{\mathcal{C}}(X,Y) = \operatorname{Ext}_{\mathcal{C}}^n(X,Y)$. Since $\operatorname{Hom}_{\mathcal{C}}$ is a bifunctor, we can choose to fix a variable and $\mathscr{R}\operatorname{Hom}_{\mathcal{C}}(X,-)(Y) \cong \mathscr{R}\operatorname{Hom}_{\mathcal{C}}(-,Y)(X)$ (Theorem 4.7.10 in \cite{Schapira}). Since the injective model structure on $\operatorname{Ch}_{\bullet}(\mathcal{C})^{\operatorname{op}}$ is actually given by projective complexes in $\mathcal{C}$, this has the advantage of letting us choose which resolutions to use to use to compute $\operatorname{Ext}$.
\end{enumerate}
\end{Examples}

\section{$G$-modules and cohomology}
\label{sec:GModules}

For starters, denoting the group ring of $G$ by $\mathbb{Z}[G]$, a (left, say) \textit{$G$-module} is a (unital left) $\mathbb{Z}[G]$-module in the usual sense of a module over a ring. This is equivalent to the data of an abelian group $M$ and a group homomorphism $G \rightarrow \operatorname{Aut_{Ab}}(M)$. A $G$-linear map $f \colon M \rightarrow N$ between $G$-modules is an abelian group homomorphism $f$ such that, for all $m\in M$ and $g \in G$, we have $f(g\cdot m) = g\cdot f(m)$. These form the category $G\operatorname{Mod}$ of $G$-modules. A $G$-module is \textit{trivial} if the action is given $g \cdot m = m$ for all $g\in G$ and $m\in M$. Some important examples of trivial $G$-modules are as follows.
\begin{Examples}\phantom{-}
\begin{enumerate}
    \item Any abelian group $A$ can be thought of as a trivial $G$-module, namely $\operatorname{triv}_G(A) = A$ with action $g \cdot a \coloneqq a$ for all $a\in A, g\in G$.
    
    \item If $M$ is any $G$-module, there is a trivial sub $G$-module by taking invariants under the $G$-action:
    \[
    M^G \coloneqq \left\{m \in M\ |\ g\cdot m = m,\ \forall g\in G \right\}
    \]

    \item If $M$ is any $G$-module, there is a trivial quotient $G$-module by taking co-invariants under the $G$-action:
    \[
    M_G \coloneqq \frac{M}{\left\langle g\cdot m - m\right\rangle \scriptstyle{_{g\in G}^{m\in M}}}
    \]
\end{enumerate}
\end{Examples}

\noindent These examples are all functorial; $\operatorname{triv}_G$ is an inclusion of categories $\operatorname{Ab}\hookrightarrow G\operatorname{Mod}$, with $(-)^G$ and $(-)_G$ being its right and left adjoints respectively. This is depicted by the following diagram.
\[\begin{tikzcd}
	{G\operatorname{Mod}} && {\operatorname{Ab}}
	\arrow[""{name=0, anchor=center, inner sep=0}, "{{(-)_G}}", curve={height=-18pt}, from=1-1, to=1-3]
	\arrow[""{name=1, anchor=center, inner sep=0}, "{{(-)^G}}"', curve={height=18pt}, from=1-1, to=1-3]
	\arrow[""{name=2, anchor=center, inner sep=0}, "{{\operatorname{triv}_G}}"{description}, hook', from=1-3, to=1-1]
	\arrow["\dashv"{anchor=center, rotate=-90}, shift right, draw=none, from=0, to=2]
	\arrow["\dashv"{anchor=center, rotate=-90}, shift right, draw=none, from=2, to=1]
\end{tikzcd}\]
\noindent By abstract nonsense, since both $\operatorname{Ab}$ and $G\operatorname{Mod}$ are abelian categories, $(-)^G$ being a right adjoint functor implies it is left exact. By the further abstract nonsense of Appendix \ref{sec:ChainHomotopy}, since it's a left exact functor, it admits a right derived functor.

\begin{Definition}
Let $G$ be a group, and $M$ some $G$-module. The \textit{$n^{th}$ cohomology group of $G$ with coefficients in $M$} is the $n^{th}$ right derived functor $\mathscr{R}^n$ of $(-)^G$.
\[
H^n(G;M) \coloneqq \mathscr{R}^n(M)
\]
Since $(M)^G \cong \operatorname{Hom}_G(\Z, M)$, where $\Z$ has the trivial $G$-module structure, then $H^n(G;M) \cong \operatorname{Ext}_{\Z[G]}^n(\Z,M)$.
\end{Definition}

As pointed out in Example \ref{ex:DerivedFuncs}.2, to compute $\operatorname{Ext}_{\Z[G]}^n(\Z,M)$, I could proceed via an injective resolution of $M$ and $\operatorname{Hom}_G(\Z,-)$, or a projective resolution of the trivial $G$-module $\Z$ and $\operatorname{Hom}_G(-,M)$. This is where I make use of the famous \textit{normalized bar resolution} $\overline{\Z G}_{\bullet}$, which is a \textit{free} resolution of $\Z$. A basis of $\overline{\Z G}_n$ is given by the symbols $[g_1|\cdots|g_n]$, for $g_i \in G\setminus \{1\}$. From this point on, any computational details may be found in any standard reference, e.g. \cite{GroupCohomology}. Here is a more operational definition of group cohomology.

\begin{Definitions}
\label{def:BarCohomology}
Let $G$ be any group.
\begin{enumerate}
    
    \item For any $n \in \mathbb{N}$, $C^n(G;M)$ is the group of \textit{normalized n-cochains of $G$ in $M$}. These are functions $f \colon G^{\times n} \rightarrow M$ such that $f(g_1,\dots,g_n) = 0$ if any $g_i = 1$, and are obtained by $\operatorname{Hom}$-ing the bar resolution into $M$.
    \begin{align*}
        \operatorname{Hom}_G(\overline{\Z G}_n, M) &\cong C^n(G;M)\\
        \varphi &\mapsto \varphi([g_1|\cdots|g_n])
    \end{align*}

    \item For any $n \in \mathbb{N}$, the resulting \textit{coboundary homomorphism} $\delta_G^n \colon C^n(G;M) \rightarrow C^{n+1}(G;M)$ is given by
    \begin{align*}
        (\delta_G^n f) (g_1,\dots,g_{n+1}) = g_1&\cdot f(g_2,\dots,g_n)+\cdots\\ 
        &+\sum_{i=1}^{n-1}(-1)^i f(g_1,\dots,g_ig_{i+1},\dots,g_n)+\cdots\\
        &+(-1)^n f(g_1,\cdots,g_{n-1})
    \end{align*}

    \item The kernel of the coboundary, $\operatorname{ker}(\delta_G^{n}) = Z^n(G;M)$, is the group of (normalized) $n$-\textit{cocycles} of $G$ in $M$. Similarly the image, $\operatorname{im}(\delta_G^{n-1}) = B^n(G;M)$, is the group of (normalized) $n$-\textit{coboundaries} $G$ in $M$.

    \item The \textit{$n^{th}$ cohomology group of $G$ with coefficients in $M$} may be computed as the quotient group
    \begin{align*}
        H^n(G;M) = \frac{Z^n(G;M)}{B^n(G;M)}
    \end{align*}
\end{enumerate}
\end{Definitions}

\noindent When the $G$-module $M$ has a non-trivial $G$-action, then this is also known as ``twisted" group cohomology. From here out, I take $G$ to be a fixed, but arbitrary, finite group and all cochains are taken to be normalized. Beware that when writing coboundary homomorphisms, I will often leave indices implicit.

\section{Pairings and Cup Products}
\label{ap:C}

Now let's focus on the coefficients. For any fixed $G$-module $M$, the cochains of $G$ in $M$ naturally form a graded abelian group, with grading given by the degree of cochains. Since cohomology is defined degree-wise, this graded structure then descends to cohomology, so $H^*(G;M) = \bigoplus_{n \geq 0} H^n(G;M)$. This graded setting allows for more structure to be found. In particular, there is the following lemma.
\begin{Lemma}
\label{lem:CupProduct}
The functor $H^*(G;-) \colon G\operatorname{Mod} \rightarrow \operatorname{Ab_{\mathbb{N}}}$ is lax monoidal. In particular, there is a natural map $\cup \colon H^*(G;M) \otimes H^*(G;N) \rightarrow H^*(G;M \otimes N)$.
\end{Lemma}
\begin{Proof}
Recall that we compute a derived functor by computing the original functor on a fibrant or cofibrant replacement. In this case, we are in the injective model structure on $\operatorname{Ch}_{\bullet}(G\operatorname{Mod})^{\text{op}}$, so in particular we have the following.
\begin{align*}
H^n(G;M) &= \operatorname{Ext}^n_{\Z[G]}(-,M)(\Z)\\
&= \operatorname{Hom}_G(-,M) \circ \mathcal{F}ib_{G\operatorname{Mod}^{\text{op}}}(\Z)\\
&= \operatorname{Hom}_G(-,M) \circ \mathcal{C}of_{G\operatorname{Mod}}^{\text{op}}(\Z)
\end{align*}
\noindent where I used that opposite of a model category gives another model category, but where the fibrations and cofibrations have swapped roles. In particlar we have $\mathcal{F}ib_{G\operatorname{Mod}^{\text{op}}}(\Z) = \mathcal{C}of_{G\operatorname{Mod}}^{\text{op}}(\Z) = \overline{\Z G}_{\bullet}$.

The trivial $G$-module $\Z$ is a trivial coalgebra in $G\operatorname{Mod}$. Specifically, $\Z \otimes_{\Z[G]} \Z = \Z \otimes_{\Z} \Z \cong \Z$, and so the identity map defines a coassociative comultiplication. Moreover, Example \ref{ex:DerivedFuncs}.1 implies that $\mathcal{C}of_{G\operatorname{Mod}}^{\text{op}}$ is oplax monoidal, and so preserves coalgebras. Indeed, the bar resolution is a similarly trivial coalgebra. This comes from the fact that the tensor product $\overline{\Z G}_{\bullet} \otimes \overline{\Z G}_{\bullet}$ is also a (projective) resolution of $\Z$, together with Theorem I.7.5 of \cite{GroupCohomology}. This latter states that any two projective resolutions over the same module unique up to a unique isomorphism in $\mathscr{D}_{\bullet}(G\operatorname{Mod})^{\text{op}}$. Finally, in general if $X \in \mathcal{C}$ is a coalgebra in a (enriched) monoidal category $\mathcal{C}$, then the functor $\operatorname{Hom}_{\mathcal{C}}(X,-)$ is lax monoidal.
\end{Proof}

In order to proceed, we need to make an explicit choice of chain map witnessing the (up to homotopy) trivial comultiplication of the bar resolution. A common choice is called the Alexander-Whitney chain map, and is defined degree-wise (e.g. in \cite[\S 5.1]{GroupCohomology}) as follows.
\begin{gather*}
    \Delta_n \colon \overline{\Z G}_n \longrightarrow \bigoplus_{i+j=n}\overline{\Z G}_i \otimes \overline{\Z G}_j\\
    [g_1|\cdots|g_n] \mapsto \sum_{k=0}^n [g_1|\cdots|g_k] \otimes g_1 \cdots g_k [g_{k+1}|\cdots|g_n]
\end{gather*}
Thus there's a natural transformation $C^{\bullet}(G;M) \otimes C^{\bullet}(G;N) \rightarrow C^{\bullet}(G;M\otimes N)$ given by $\alpha \otimes \beta \mapsto \alpha \cup \beta \coloneqq (\alpha \otimes \beta) \circ \Delta_{\bullet}$. Moreover, this map satisfies the \textit{Leibniz rule}, given by the following equation.
\begin{equation}
\label{eq:Leibniz}
    \delta (\alpha \cup \beta) = (\delta \alpha) \cup \beta + (-1)^p \alpha \cup (\delta \beta)
\end{equation}
In particular, this means that $\cup$ sends cocycles to cocylces and coboundaries to coboundaries, and so is well-defined in cohomology. This map has many nice properties that are collected in \cite[\S 5.3]{GroupCohomology}, including being graded-commutative, associative, unital. The formula in (\ref{eq:Leibniz}) is not hard to prove in general, but the notation gets extremely cumbersome. Here is a sample computation with $\alpha, \beta \in C^{1}(G;M)$.
\begin{align*}
    \delta (\alpha \cup \beta) (g,h,k) &= g\cdot(\alpha(h)\otimes h\cdot\beta(k)) - \alpha(gh) \otimes (gh)\cdot\beta(k)\\
    &\ \ \ \ + \alpha(g)\otimes g \beta(hk) - \alpha(g) \otimes g\cdot \beta(h)\\
    &= g\cdot\alpha(h)\otimes (gh)\cdot\beta(k) - \alpha(gh) \otimes (gh)\cdot\beta(k)\\
    &\ \ \ \ + \alpha(g)\otimes g \beta(hk) - \alpha(g) \otimes g\cdot \beta(h)\\
    &\ \ \ \ + \textcolor{black}{\alpha(g) \otimes (gh)\cdot\beta(k)} - \textcolor{black}{\alpha(g) \otimes (gh)\cdot\beta(k)}\\
    &= \delta\alpha(g,h)\otimes(gh)\cdot \beta(k) + \textcolor{black}{(-1)^1} \alpha(g) \otimes g\cdot \delta\beta(h,k)\\
    &= \Big((\delta \alpha) \cup \beta + (-1)^1 \alpha \cup (\delta\beta)\Big)(g,h,k)
\end{align*}

We can use the natural map implied by Lemma \ref{lem:CupProduct} to find some more structure in graded group cohomology. In particular, consider three $G$-modules $M,N,L$, and let $G$ act diagonally on $M\otimes N$ i.e. $g\cdot (m\otimes n) = (g\cdot m) \otimes (g\cdot n)$. This situation lets us consider some particularly nice graded homomorphisms of the form $H^*(G;M) \otimes H^*(G;M) \rightarrow H^*(G;L)$.

\begin{Definitions}
\label{def:CupProduct}
A \textit{pairing} of $M$ and $N$ to $L$ is a $G$-module homomorphism $P \colon M\otimes N \rightarrow L$. The \textit{cup product} associated to the pairing $P$ is the composition $\cup_P = H^*(P) \circ \cup$.
\end{Definitions}

\noindent In section \ref{sec:SpecSeq}, I recall a simple, yet important, example of a pairing whose associated cup product plays in important role in the Lyndon-Hochschild-Serre spectral sequence. Here are some other simple examples of cup products.
\begin{Examples}\phantom{-}
    \begin{enumerate}
        \item Trivially, taking $L = M\otimes N$, the pairing $\operatorname{id}_{M\otimes N}$ induces $\cup_{\operatorname{id}} = \cup$.

        \item Let $M=N=L=R$, some ring $R$, and let $G$ act trivially. The map $(r,s) \mapsto rs$ is a pairing, and the resultant cup product gives $H^*(G;M)$ a graded ring structure. This is exactly what provides the cup product in the cohomology ring of a space in algebraic topology, hence the name in the more general setting.
    \end{enumerate}
\end{Examples}

\sloppy
\printbibliography

\end{document}